\documentclass[
    ,final            % use final for the camera ready runs
  ]
  {aipproc}

\layoutstyle{8x11single}

\begin{document}

\title{A semi-implicit high-order space-time scheme on staggered meshes for the 2D incompressible Navier-Stokes equations}


\author{Francesco Lohengrin Romeo}{
  address={Laboratory of Applied Mathematics,\\Department of Civil, Environmental and Mechanical Engineering,\\
  	University of Trento, \\
  	Via Mesiano, 77, I-38123 Trento, Italy }
}

%\author{<author2>}{
%  address={<common address for author2 and author3>}
%}

%\author{<author3>}{
%  address={<common address for author2 and author3>}
%  ,altaddress={<author1 address>} % additional visiting address
%}

\begin{abstract}
A new high order accurate semi-implicit space-time Discontinuous Galerkin method on staggered grids, for the simulation of viscous incompressible flows on two-dimensional domains is presented. The designed scheme is of the Arbitrary Lagrangian Eulerian type, which is suitable to work on fixed as well as on moving meshes. In our space-time formulation, by expressing the numerical solution in terms of piecewise space-time polynomials, an arbitrary high order of accuracy in time is achieved through a simple and efficient method of Picard iterations. For the dual mesh, the basis functions consist in the union of continuous piecewise polynomials on the two subtriangles within the quadrilaterals: this allows the construction of a quadrature-free scheme, resulting in a very efficient algorithm. Some numerical examples confirm that the proposed method outperform existing ones.
\end{abstract}

\maketitle

%%%%%%%%%%%%%%%%%%%%%%%%%%%%%%%%%%%%%%%%%%%%
%% MAINMATTER
%%%%%%%%%%%%%%%%%%%%%%%%%%%%%%%%%%%%%%%%%%%%

\section{The new staggered space-time finite element method}

\subsection{The spatial discretization}
We propose a numerical scheme for the solution of the incompressible Navier-Stokes equations on 2D domains, following the original ideas in \cite{2DSIUSW,2STINS}. The spatial discretization is performed through the use of two unstructured grids: a \textit{primary} triangular grid for the approximated pressure function, and a staggered edge-based quadrilateral grid (named \textit{dual} grid) for the approximated velocity functions. In the framework of the Galerkin finite element methods, we define the finite element space of order $p$ for the discretized pressure using the standard nodal basis functions on the reference triangular element $T_{std}=\{(\xi,\eta) \in R^{2,+} \,\, | \,\, \eta\leq1-\xi \wedge 0 \leq \xi \leq 1 \}$, by imposing the classical Lagrange interpolation condition $\phi_k(\vec \xi_l) = \delta_{kl}$ over the 2D Newton-Cotes quadrature nodes:
\begin{equation}
\vec \xi_k = (\xi_{k_1}, \gamma_{k_2}) = \left( \frac{k_1}{p}, \frac{k_2}{p} \right),
\label{eqn.nodes}
\end{equation}
with the multi-index $k=(k_1,k_2)$ and the index ranges $0 \leq k_1 \leq p$ and $0 \leq k_2 \leq p - k_1$. In this way $N_\phi=\frac{(p+1)(p+2)}{2}$ basis
functions $\{\phi_k \}_{k \in [1,N_\phi]}$ are obtained.
Analogously, we choose the $N_\psi=(p+1)^2$ basis functions $\{\psi_k\}_{k \in [1, N_\psi]}$ on the reference square element $R_{std}=\{(\xi,\eta) \in R^{2,+} \,\, | \,\, 0 \leq \xi \leq 1 \wedge 0 \leq \eta \leq 1 \}$ for the dual finite element space of order $p$. Since we want to derive a quadrature-free Arbitrary Lagrangian Eulerian implementation, we consider the square as the union of two sub-triangles $ T_{I} $ and $ T_{II} $ and we construct the basis functions following the standard nodal approach of continuous finite elements (indeed, we get two mini-elements). This is an approach very similar to the $ {P}^{k}-iso$ $ {P}^{k+1} $ finite elements used for the velocity approximation in some mixed problems, for example the Stokes problem, where mixed finite element approximations are employed in order to numerically satisfy the inf-sup compatibility condition. These elements are fully described in \cite{ern2004}.
We build the basis functions over the nodes which lie on the sub-triangle $ T_{I} $ exactly like we did for $ T_{std} $, and then we extend them continuously to $ T_{II} $, in order not to have discontinuity inside the square. Viceversa, we construct the basis functions over the nodes which lie on $ T_{II} $ via a transformation between $ T_{II} $ and $ T_{std} $, and then we extend them continuously to $ T_{I} $.

\subsection{Space-time extension}
For the time discretization, the generic space-time element defined in the $n$-th time-step of the simulation $[t^n, t^{n+1}]$ is given by a triangular prism for the main grid, and a quadrilateral prism for the dual grid. In the Lagrangian case, these volumes are stretched in some way according to the local mesh velocity, while in the Eulerian case they are some right, not slanting in time, prisms.
The $N_\gamma=p_\gamma+1$ temporal basis functions $\{\gamma_k\}_{k \in [1, N_\gamma]}$ for polynomials of degree $p_\gamma$ are defined as the Lagrange interpolation polynomials passing through the equidistant 1D Newton-Cotes quadrature nodes on the reference interval $I_{std} = [0,1]$. Finally, using the tensor product, we define the basis functions on the space-time elements as $\tilde{\phi}(\xi,\eta,\tau)=\phi(\xi, \eta) \cdot \gamma(\tau)$ and $\tilde{\psi}(\xi, \eta,\tau)=\psi(\xi, \eta) \cdot \gamma(\tau)$. Therefore, the total number of basis functions becomes $N_\phi^{st}=N_\phi \cdot N_\gamma$ and $N_\psi^{st}=N_\psi \cdot N_\gamma$.

\section{The semi-implicit scheme for the incompressible Navier-Stokes equations}

We consider the two dimensional Navier-Stokes equations for incompressible Newtonian fluids with homogeneous density $\rho$ and homogeneous dynamic viscosity coefficient $\mu$, in the conservative, adimensional form:
\begin{eqnarray}
\nabla \cdot \mathbf{v} & = & 0 \label{eq:CS_2}, \\
\frac{\partial \mathbf{v}}{\partial t}+\nabla \cdot \mathbf{F}_c + \nabla p & = & \nu \Delta \mathbf{v} + \mathbf{S} \label{eq:CS_2_2_0}, 
\end{eqnarray}
where $p=P/\rho$ indicates the normalized fluid pressure; $P$ is the physical pressure; $\nu = \mu / \rho$ is the kinematic viscosity coefficient; $\mathbf{v}=\left( u , v \right)^T$ is the velocity vector; $u$ and $v$ are the velocity components in the $x$ and $y$ direction, respectively;
$\mathbf{S}=\mathbf{S}(\mathbf{v},x,y,t)$ is a source term; 
$\mathbf{F}_c=\mathbf{v} \otimes \mathbf{v}$ is the flux tensor of the nonlinear convective terms.
Following the same ideas in 
\cite{ADERNSE,MunzDiffusionFlux}, the momentum Eq. \eqref{eq:CS_2_2_0} can be rewritten as:
\begin{equation}
\frac{\partial \mathbf{v}}{\partial t}+\nabla \cdot \mathbf{F} + \nabla p= \mathbf{S} 
\label{eq:CS_2_2},
\end{equation}
where $\mathbf{F}=\mathbf{F}(\mathbf{v},\nabla \mathbf{v})=\mathbf{F}_c(\mathbf{v})-\nu \nabla \mathbf{v}$.

\subsection{The Picard's method for the high-order accuracy in time}

The Discontinuous Galerkin finite element formulation considers the integration of Equations \eqref{eq:CS_2} and \eqref{eq:CS_2_2_0} on the space-time control volumes of the primary grid and of the dual grid, respectively. Following the ideas in \cite{casulli}, \cite{2STINS}, a \textit{semi-implicit} discretization is employed, which combines
the simplicity of explicit methods for nonlinear hyperbolic PDE with the stability and efficiency of implicit time discretizations.
In order to obtain a method with high-order accuracy in time, we use a simple \textit{Picard iteration} which introduces the information of the new pressure into the viscous and convective terms, but without involving a nonlinearity in the final system to be solved. This approach is inspired by the local space-time Galerkin predictor method proposed for the high order time discretization of $P_NP_M$ schemes in \cite{Dumbser2008,ADERNSE}. One time step of the final numerical scheme can be summarized as follows: 
\begin{enumerate}
	\item Initialize $\mathbf{v}_h^{n+1,0}$ and $p_h^{{n+1,0}}$ using the known information from the previous time-step, or the initial conditions; 
	\item Picard iteration over $k=0, \ldots, N_{pic}$:
	\begin{enumerate}
		
		\item compute $\mathbf{v}_h^{n+1,k+1/2}$ using $p_h^{n+1,k}$ in the discretized momentum equation; then 
		set $\mathbf{Fv}_h^{n+1,k +1/2}:=\mathbf{v}_h^{n+1,k+1/2}$;
		 
		\item compute $p_h^{{n+1,k+1}}$ by solving the discrete pressure Poisson equation, coming from a formal substitution of the velocity unknowns from the momentum equation \eqref{eq:CS_2_2_0} into the incompressibility constraint \eqref{eq:CS_2};  
		\item update $\mathbf{v}_{h}^{n+1,k+1}$ from the pressure correction values $\Delta p_h^{n+1,k+1}=p_h^{n+1,k+1}-p_h^{n+1,k} $;  
		
	\end{enumerate}	
	\item set $\mathbf{v}_{h}^{n+1}=\mathbf{v}_{h}^{n+1,k+1}$	and  $p_{h}^{{n+1}}=p_{h}^{{n+1,k+1}}$. 
\end{enumerate}
The resulting linear system for the pressure correction is very sparse thanks to the use
of the \textit{staggered} grid, including only four non-zero blocks per element.

\section{Numerical results}
Some relevant tests were executed in order to assess the computational efficiency and the accuracy of the new numerical method. 
Compared to the staggered space-time DG algorithm of Tavelli and Dumbser \cite{2STINS} the new method proposed here is not only computationally more efficient thanks to its quadrature-free formulation, but also less memory consuming, since all integrals can be precomputed once and for all on a universal reference element. 
Moreover, thanks to the use of the piecewise basis functions for the dual finite element space, all the matrices of the Galerkin formulation can be updated, with their geometric information, at every time-step $n$, by a cheap matrix-vector product which uses the first levels of the cache memory.

\begin{figure}
  \includegraphics[height=.23\textheight]{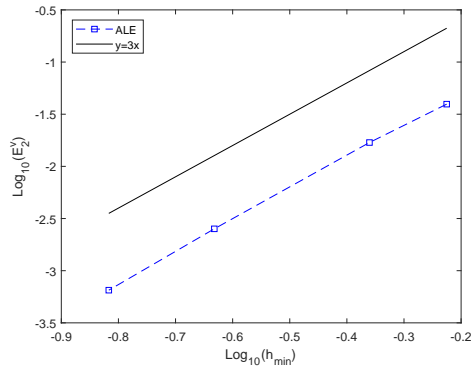}
  \caption{Taylor-Green vortex test: slope of the third-order method (convergence order for the velocity).}
  \label{TG_P2_P2_conv}
\end{figure}

\begin{figure}
	\includegraphics[height=.23\textheight]{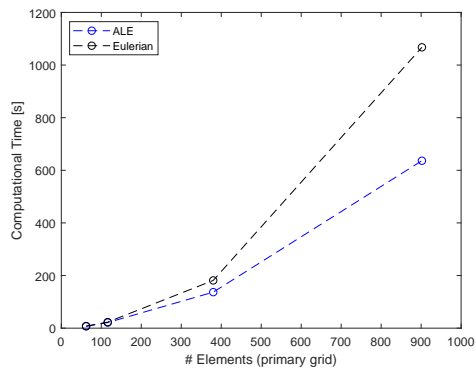}
	\caption{Taylor-Green vortex test: computational times required by the third-order method.}
	\label{TG_P2_P2_CPUtimes}
\end{figure}

\paragraph{Taylor-Green vortex test}
\begin{table}
	\begin{tabular}{lrrrrrr}
		\hline
		& \tablehead{1}{r}{b}{$E_{2}^{p}$}
		& \tablehead{1}{r}{b}{$E_{2}^{p}$ [Eul.]}
		& \tablehead{1}{r}{b}{$E_{2}^{v}$}
		& \tablehead{1}{r}{b}{$E_{2}^{v}$ [Eul.]}
		& \tablehead{1}{r}{b}{$\sigma_{2}^{v}$}
		& \tablehead{1}{r}{b}{$\sigma_{2}^{v}$ [Eul.]}   \\
		\hline
		62   & 6.26E-02 & 4.40E-02 & 3.95E-02 & 2.18E-02 & - & - \\
		116   & 2.67E-02 & 1.87E-02 & 1.69E-02 & 8.43E-03 & 2.7 & 3.1 \\
		380   & 4.49E-03 & 3.57E-03 & 2.52E-03 & 1.17E-03 & 3.0 & 3.2 \\
		902   & 1.36E-03 & 1.11E-03 & 6.50E-04 & 2.97E-04 & 3.2 & 3.2 \\
		\hline
		\caption{P2P2}
	\end{tabular}
	\caption{Taylor-Green vortex test, $p = p_{\gamma} =2$: errors and convergence order.}
	\label{tab:a2}
\end{table}
\begin{table}
	\begin{tabular}{lrrrrrr}
		\hline
		& \tablehead{1}{r}{b}{CPU time}
		& \tablehead{1}{r}{b}{CPU time [Eul.]}
		& \tablehead{1}{r}{b}{RAM usage}
		& \tablehead{1}{r}{b}{RAM usage [Eul.]}  \\
		\hline
		62   & 7 & 7 & 49 & 67\\
		116   & 21 & 22 & 67 & 106 \\
		380   & 137 & 178 & 179 & 311\\
		902   & 636 & 1059 & 392 & 709 \\
		\hline
		\caption{P2P2}
	\end{tabular}
	\caption{Taylor-Green vortex test, $p = p_{\gamma} =2$: CPU times (in seconds) and memory consumption (in MB).}
	\label{TG_P2_P2_eff}
\end{table}
As a numerical test for the incompressible Navier-Stokes equations, we consider the unsteady Taylor-Green vortex problem. The solution of the problem is determined by:
\begin{eqnarray}
p(x,y,t) =  \frac{1}{4} \left[\cos(2x)+\cos(2y) \right] \exp(-4\nu t) \label{eq:TG_p}, \\
u(x,y,t) = \sin(x)\cos(y)\exp(-2\nu t) \label{eq:TG_u}, \\
v(x,y,t) = -\cos(x)\sin(y)\exp(-2\nu t) \label{eq:TG_v}, 
\end{eqnarray}
We set periodic boundary conditions for the domain $ \Omega = [-\pi,\pi]^2 $, $\nu=0.1$ and the final time $T = 0.1$. An implicit treatment of the viscosity terms is applied, therefore the time step is given by the CFL-type restriction for only the convective operator:
$$
\Delta t = \frac{CFL}{2p+1} \cdot \frac{h_{min}}{s_{max}},
$$
where $CFL=0.4$, $h_{min}$ is the minimum of the radii of the circles inscribed in the triangles (primary grid) and $s_{max}$ is the maximum, over all the edges of the quadrilaterals (dual grid), of the maximum eigenvalue of the convective flux tensor, that is $2 \  max\{|\mathbf{v}^{+} \cdot \mathbf{n}|,|\mathbf{v}^{-} \cdot \mathbf{n}|\}$.\\
We have compared the results obtained by the new ALE implementation of the method, using zero velocity mesh, with the results that were obtained in \cite{2STINS} by an Eulerian implementation of the method, with fixed meshes: in Table \ref{tab:a2} the errors in the $ L^2 $ norms for the pressure $p$ and the velocity field $\mathbf{v}$, together with the convergence order $\sigma_{2}^{v} = \frac{\log\left( (E_{2}^{v})_{1}/(E_{2}^{v})_{2} \right)}{\log\left((h_{min})_{1} / (h_{min})_{2} \right)}$, for an increasing size of the primary grid (first column), are reported for the case $p=p_{\gamma}=2$.
In Figure \ref{TG_P2_P2_conv} the slope of the third-order method is shown for this test.\\ Finally, Table \ref{TG_P2_P2_eff} and Figure \ref{TG_P2_P2_CPUtimes} show the improved efficiency of the new algorithm (ALE, with zero mesh velocity), with respect to the Eulerian implementation: for the grid with the highest resolution, $39.9\%$ of the computational time is saved, and only $44.7\%$ of the memory is required, with respect to the Eulerian implementation. 

%%%%%%%%%%%%%%%%%%%%%%%%%%%%%%%%%%%%%%%%%%%%%%%%
%% BACKMATTER
%%%%%%%%%%%%%%%%%%%%%%%%%%%%%%%%%%%%%%%%%%%%%%%%

\begin{theacknowledgments}
  This work was partially funded by the research project \textsl{STiMulUs}, ERC Grant agreement no. 278267. Financial support has also been provided by the Italian Ministry of Education, University and Research (MIUR) in the frame of the \textit{Departments of Excellence Initiative 2018-2022} attributed to DICAM of the University of Trento (grant L. 232/2016). The author has also received funding from the University of Trento via the Strategic Initiative \textit{Modeling and Simulation}.
\end{theacknowledgments}



\bibliographystyle{aipproc}   % if natbib is available

%\bibliography{romeo}

\IfFileExists{\jobname.bbl}{}
 {\typeout{}
  \typeout{******************************************}
  \typeout{** Please run "bibtex \jobname" to optain}
  \typeout{** the bibliography and then re-run LaTeX}
  \typeout{** twice to fix the references!}
  \typeout{******************************************}
  \typeout{}
 }

%%%%%%%%%%%%%%%%%%%%%%%%%%%%%%%%%%%%%%%%%%%
%% The following lines show an example how to produce a bibliography
%% without the help of the BibTeX program. This could be used instead
%% of the above.
%%%%%%%%%%%%%%%%%%%%%%%%%%%%%%%%%%%%%%%%%%%

\end{document}